\documentclass{amsart}
\usepackage{graphicx}
\usepackage{amssymb}
\usepackage{epstopdf}
\usepackage{nicefrac}
\usepackage{enumerate}
\usepackage{float}
\usepackage{color}
\usepackage{pst-all}
\usepackage{tabularx}
\usepackage{hyperref}
\usepackage{caption}
\usepackage{subcaption}

\newtheorem{thm}{THEOREM}[section]

\newtheorem{prob}[thm]{PROBLEM}

\newtheorem*{prob*}{PROBLEM}

\newcommand{\e}{{\epsilon}}

\newcommand{\mA}{{\mathbb A}}

\newcommand{\mK}{{\mathbb K}}

\newcommand{\mS}{{\mathbb S}}
\newcommand{\mT}{{\mathbb T}}

\newcommand{\mW}{{\mathbb W}}

\newcommand{\cC}{{\mathcal C}}

\newcommand{\cK}{{\mathcal K}}

\newcommand{\cO}{{\mathcal O}}

\newcommand{\cR}{{\mathcal R}}

\newcommand{\cW}{{\mathcal W}}

\newcommand{\fD}{{\mathfrak{D}}}
\newcommand{\fM}{{\mathfrak{M}}}

\begin{document}

\title{Kuperberg dreams}

\author{Steven Hurder}
\address{Steven Hurder, Department of Mathematics, University of Illinois at Chicago, 322 SEO (m/c 249), 851 S. Morgan Street, Chicago, IL 60607-7045}
\email{hurder@uic.edu}
\thanks{Preprint date: December 8, 2021; to appear online at  Celebration.org}

\author{Ana Rechtman}
\address{Ana Rechtman, Institut de Recherche Math\'{e}matique Avanc\'{e}e,   Universit\'{e} de Strasbourg, 7 rue Ren\'{e} Descartes, 67084 Strasbourg, France}
\email{rechtman@math.unistra.fr}

\date{}



 \begin{abstract}
The ``Seifert Conjecture'' stated, ``Every non-singular vector field on the 3-sphere ${\mathbb S}^3$ has a periodic orbit''.  
In a celebrated work, Krystyna Kuperberg   gave a construction of a smooth aperiodic vector field on a plug, which is then used to construct    counter-examples to the  Seifert Conjecture for smooth flows on the 3-sphere, and on compact 3-manifolds in general. The dynamics of the flows in these plugs have been extensively studied, with more precise results known  in special ``generic'' cases of the construction. Moreover, the dynamical properties of smooth perturbations of Kuperberg's construction have been considered.  
    In this work, we discuss the genesis of Kuperberg's construction as an evolution from the Schweitzer construction of an aperiodic plug. We discuss some of  the known results for Kuperberg flows, and   discuss some of the many interesting questions and problems that remain open about their dynamical and ergodic properties.
 \end{abstract}

  \maketitle

 \section{Birth of a method}\label{sec-intro}
In August   1993, Krystyna Kuperberg wrote out a  3 page sketch, of a construction of  what she thought would be a $C^3$-counterexample to the Seifert Conjecture. The idea for this construction had occurred to her earlier in the month, while attending the \emph{Georgia Topology Conference} in Athens, Georgia. She faxed these notes to her son Greg, with the admonition that he  wait for a more final draft.
In the following days, she   wrote out the precise construction, which expanded to an 8 page note, and realized that it was a $C^{\infty}$-counterexample.   
Greg then   sent an email to several leading topologists, announcing the exciting news, that his mother had found a smooth counterexample to the Seifert Conjecture!

The announcement of Krystyna's work caused a sensation. William Thurston  lectured on her construction at the MSRI, Berkeley in September, and observed that the construction could be realized as a \emph{real analytic} aperiodic flow. Krystyna   mailed around 60 copies of the manuscript to mathematicians around the world. Then  in October, she gave a plenary lecture on her work at the conference in honor of Morris Hirsch in Berkeley California.
 Elsewhere, a seminar was held in Tokyo in November, as part of the   \emph{International Symposium/Workshop: Geometric Study of Foliations}. The participants in this informal seminar included, among many others, Shigenori Matsumoto,  \'{E}tienne Ghys, Paul Schweitzer,  and the first author.

 The manuscript by 
  Kuperberg  was   published almost immediately in the \textbf{Annals of Mathematics} \cite{Kuperberg1994}.  Moreover, soon afterwards, Ghys     presented Kuperberg's work in  a S{\'e}minaire Bourbaki \cite{Ghys1995} in June 1994.  Matsumoto   wrote a report in Japanese   \cite{Matsumoto1995}, also in 1994, which provided   more results about the dynamical properties of these ``Kuperberg flows'', with detailed proofs of their properties. 
The joint paper    by Greg and Krystyna    \cite{Kuperbergs1996} contained many new ideas and variations on the basic construction. The brief note \cite{Kuperberg1999} gives an overview of Krystyna's work on aperiodic flows. Finally,   Krystyna  Kuperberg  \cite{Kuperberg1998} gave a report to the \emph{International Congress of Mathematicians} in 1998.

  Krystyna   commented on her work: ``once you see the construction, you see it''. In addition, her original sketches of the construction, involving a thickened annulus  with ``rabbit ears'', were replaced in the published version by the now very familiar images drawn by her husband, W\l\/odzimierz (W{\l}odek), which conveyed a strong intuitive feel for the construction.  Krystyna Kuperberg was a  student of Karol Borsuk in Warsaw, and many of her works correspondingly use a strong geometric approach to analyzing problems in dynamical systems, and this is more than true for her construction of the Seifert counter-examples.

 The goals of this essay are to discuss the antecedents to  the Kuperberg construction, and give some inkling of the ideas that led Krystyna to her construction. In addition, we 
 speculate on extensions of the construction and its future, especially in the context  of other aspects of the theory of smooth dynamical systems. First, in Section~\ref{sec-past} we recall the radical idea that made the Kuperberg construction possible, as it was a break from the past approaches to constructing counterexamples in a crucial manner. Section~\ref{sec-dynamics} gives a concise overview of   known results about the Kuperberg flow.  We then discuss in Section~\ref{sec-problems}   questions about   alternate  flow dynamics for non-generic constructions of   Kuperberg flows.     Section~\ref{sec-speculations} discusses speculations and open problems for flows that preserve some additional geometric structure.

  \section{Shape of Plugs}\label{sec-past}

 In his  1950 work  \cite{Seifert1950},  Seifert introduced    an invariant for  deformations of non-singular
  flows with a closed orbit  on a $3$-manifold,    which he used to show that every sufficiently small deformation of the Hopf flow on the $3$-sphere $\mS^3$ must have a closed orbit. He also remarked in    Section~5 of this work:
\begin{quotation}
\emph{It is unknown if every continuous (non-singular) vector field on the three-dimensional sphere contains a closed integral curve.}
\end{quotation}
This remark became the basis for what is known as the  ``Seifert Conjecture'': 
\begin{quotation}
\emph{Every non-singular vector field on the $3$-sphere $\mS^3$ has a periodic orbit.}
\end{quotation}

There must exist at least one minimal set for a flow on a compact manifold - a closed invariant subset  which is minimal for this property.  A periodic orbit is a minimal set, but if there are no periodic orbits for a flow, then a minimal set for the flow will have a possibly far more    complicated topology. So one can ask if  there are restrictions on   the \emph{shape} of a minimal set for an aperiodic flow? 
 Karol Borsuk defined the shape of a topological space as an inverse limit continuum defined by  a sequence of ANR approximations of the space \cite{Borsuk1968,Borsuk1975,MardesicSegal1982,Mardesic1999}.  
For an aperiodic flow on a 3-manifold, it seems that the  minimal set  must have a complicated shape, though just how complicated is   an open question.
 
  Kuperberg's    strategy for the construction of an aperiodic flow  was not based on the choice of an invariant minimal set with prescribed shape, but was instead to create a dynamical framework which avoided the implications of the Brouwer Fixed-Point Theorem.  We explain this remark in the following discussions.

  The first general approach to the Seifert Conjecture was made in the 1966 paper by 
Wesley Wilson  \cite{Wilson1966}, where he showed that every closed $3$-manifold $M$ has a flow with only finitely many closed orbits.   Wilson introduced    the use  of a special kind of ``plug'' to modify a given flow, and this idea  profoundly influenced the thinking about the constructions of flows with special properties, and the Seifert conjecture in particular. 

\begin{figure}[!htbp]
\centering
\begin{subfigure}[c]{0.4\textwidth}{\includegraphics[height=40mm]{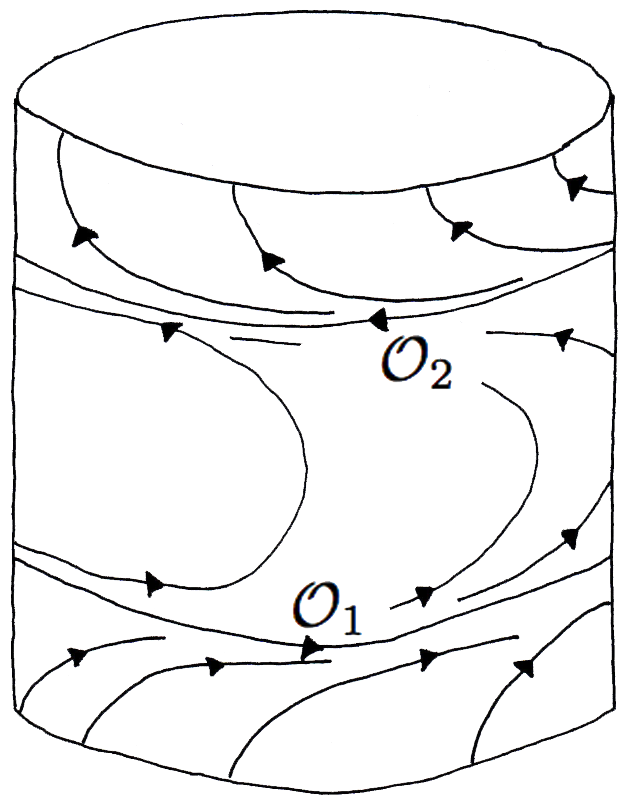}}\end{subfigure}
\begin{subfigure}[c]{0.4\textwidth}{\includegraphics[height=40mm]{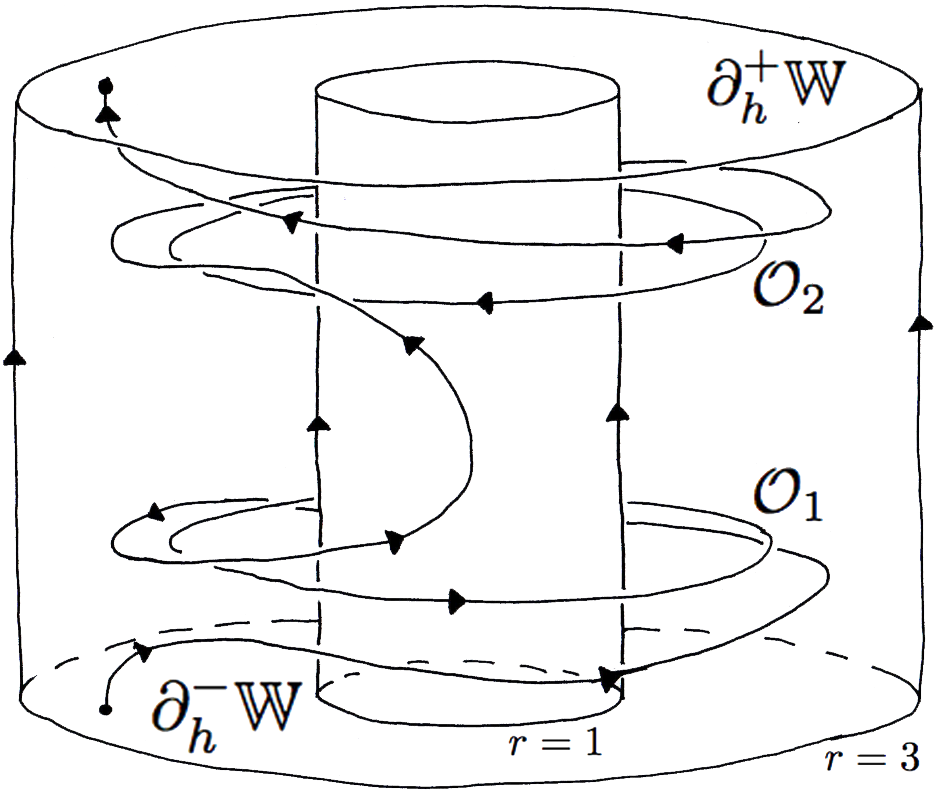}}\end{subfigure}
 \caption[{$\cW$-orbits in the cylinder  $\cC=\{r=2\}$}]{ An example of a Wilson plug with two periodic orbits \label{fig:Reebcyl}}  
\vspace{-6pt}
\end{figure}

  A \emph{Wilson plug}   for a  flow on the 3-manifold $M$  is a neighborhood $P \subset M$   diffeomorphic to  {\bf an annulus} cross an interval, $P \cong \mA^2 \times [-1,1]$,  where orbits enter one end of the plug, $\partial^- P = \mA^2 \times \{-1\}$,  and exit the other end, $\partial^+P = \mA^2 \times \{+1\}$.   The key idea introduced by Wilson was that the plug flow should be symmetric: an orbit which traverses the plug should enter and exit at symmetric points on the faces. No such restriction is placed on an orbit which enters the plug and never exits, of course, and one assumes there are such orbits. Wilson introduced the technique of using a plug reflected on itself to achieve this symmetry property. As a result, the plugs he uses have two periodic orbits, as illustrated in Figure~\ref{fig:Reebcyl}.
  
 A plug      can be inserted into any non-singular flow on a manifold, not just once but any finitely number of times as needed,   to modify a given non-singular flow.  Wilson showed that any flow can be modified  using finitely many plugs    to obtain a flow with only isolated periodic orbits.   It is an observation in \cite{PercelWilson1977} that by appropriately arranging the insertion of Wilson plugs, one can obtain a flow with exactly two periodic orbits.

\begin{figure}[!htbp]
\centering
 \includegraphics[width=50mm]{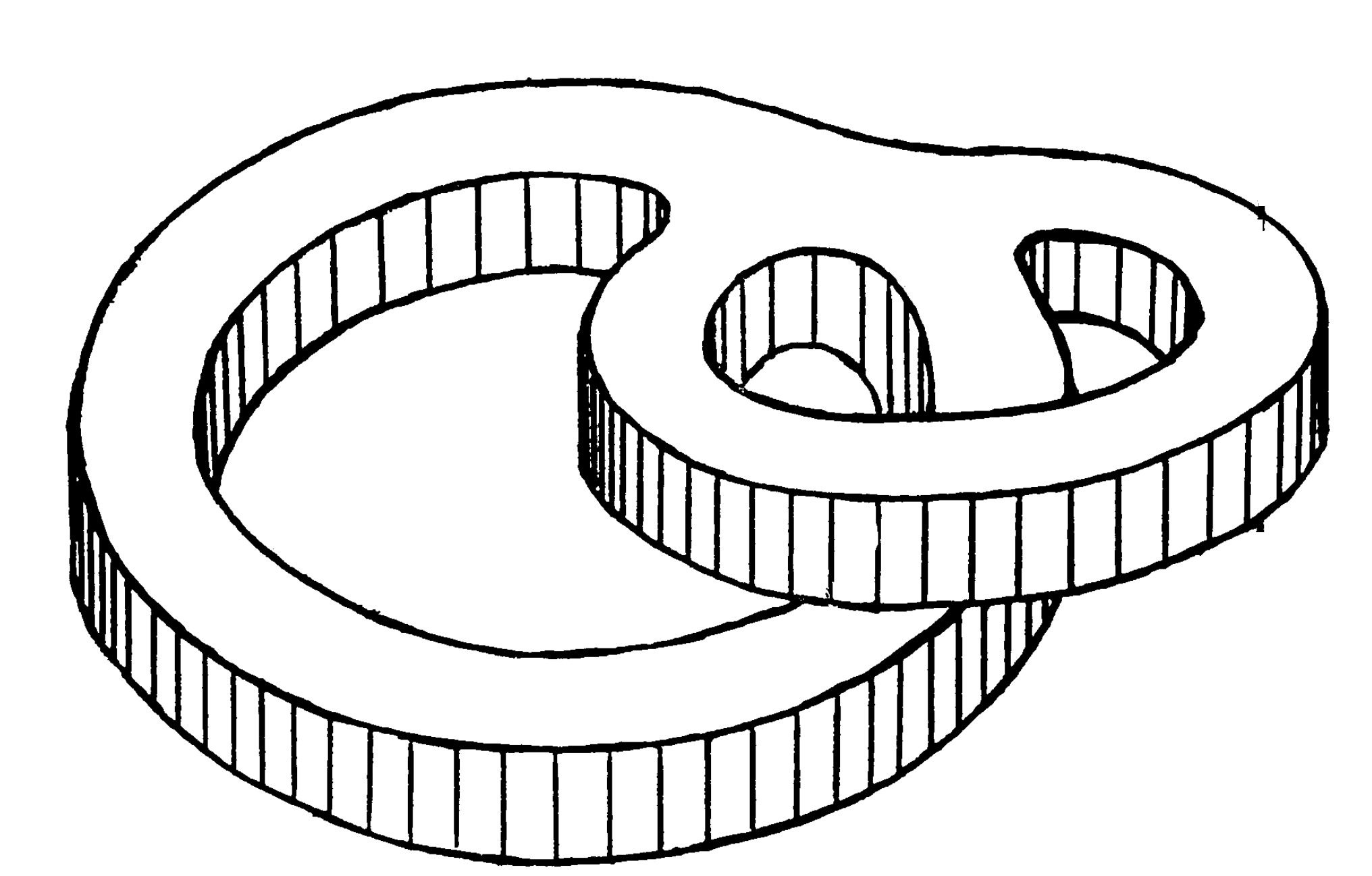}
\caption{The Schweitzer plug \label{fig:schweitzer}}
\vspace{-6pt}
\end{figure}

Subsequently, in the early 1970's, Paul  Schweitzer  had the  inspired idea to construct a modified Wilson-type plug, but instead of the top and bottom boundaries of the plug being annuli, they were diffeomorphic to a 2-torus with an open disk removed, as    illustrated in Figure~\ref{fig:schweitzer}. The paradigm shift of a Schweitzer plug is that the faces of the plug are immersed in the transverse space, but still transverse to the flow. Using a plug with this shape, the second key point was to  double the construction using a mirror symmetry as with the Wilson plug, so that   the flow in the plug has  two Denjoy-type minimal sets, instead of two circles.

A Denjoy minimal set, as introduced by Denjoy in \cite{Denjoy1932}, has the shape of a wedge of two circles.  The Schweitzer plug has  the celebrated ``double doughnut'' shape as seen in Figure~\ref{fig:schweitzer}. 
This plug  can then be inserted into flows on 3-manifolds to create new flows which have all minimal sets of Denjoy type. In particular, using this plug one can construct flows with  no periodic orbits.  Thus, Schweitzer showed that every 3-manifold carries a $C^1$-flow with no periodic orbits! This celebrated  result appeared  in the \textbf{Annals of Mathematics}   \cite{Schweitzer1974}, and  in his  talk  \cite{Schweitzer1975} to the \emph{International Congress of Mathematicians} in 1974.

  Mike Handel showed in his 1980 paper in the \textbf{Annals of Mathematics}  \cite{Handel1980} a result which severely limits the shape of a minimal set for a smooth aperiodic flow -- it cannot be embedded as the minimal set for a flow on a surface.  Thus, a flow with a   ``Denjoy type'' minimal set carried by a surface cannot be $C^2$. This was an attempt to prove that plugs were not the method to find smooth aperiodic flows on $\mS^3$.

Jenny Harrison showed in   \cite{Harrison1988} that the     Schweitzer plug could be modified to obtain  an aperiodic  plug with   a $C^{2+\delta}$-regular flow, by embedding the Denjoy minimal sets in a   ``3-dimensional''  fashion.  However, this result required delicate arguments, and obtaining any further improvement in the differentiability of the flow using  this method, seemed completely out of reach.
The  Handel and Harrison works suggested that any counter-example to the Seifert Conjecture was going to require a radically new approach.

Meanwhile, in 1979, the Kuperbergs  were attending the  \emph{Scottish Book  Conference}  in Denton, Texas, organized by Dan Mauldin  \cite{Mauldin2015}.  At this conference,  Stanislav Ulam  told Krystyna and Wlodek  about  an unsolved   problem, posed by   Ulam:
\begin{quotation}\cite[Problem~110]{Mauldin2015} 
 \emph{Let $M$ be a manifold. Does there exist a numerical constant $K$ such that every continuous mapping $f \colon M \to M$    which satisfies the condition $|f^n(x) -x| < K$  for $n=1,2,\cdots$  must    possess  a fixed point: $f^n(x_0) = x_0$?
}
\end{quotation}

The mathematician/inventor Coke Reed was   a colleague of Ulam at Los Alamos, and later a colleague of Krystyna at Auburn University.   Kuperberg and Reed discussed this problem, and subsequently solved it. They gave  counter-examples to the assertion in  Problem~110,  as published in 
the  1981 work \cite{KuperbergReed1981}. The  methods they developed are analogous to the methods introduced by  Wilson in \cite{Wilson1966}. 
 Kuperberg and Reed   improved their result in the 1989 paper \cite{KuperbergReed1989} via an application of the Schweitzer plug. An analysis of the solutions to the Ulam Problem~110 were   recalled and analyzed in the paper  
\cite{KKPR1993} by 
{K.~Kuperberg, W.~Kuperberg, P.~Minc and C.~Reed}. 
There is a connection between the solutions of the Ulam Problem and the construction of counter-examples to the Seifert Conjecture, in that the plug technique can be applied to construct counter-examples to both, as remarked  in the joint work of Greg and Krystyna \cite[Theorem~8]{Kuperbergs1996}.

Following the successful solution to the Ulam Problem 110, Kuperberg began to consider what was required to construct a smooth plug without periodic points for the flow. As she remarked \cite{Kuperberg2021}, the problem was always  the \emph{Brouwer Fixed-Point Theorem}: any map of a (transverse) disk to itself must have a fixed point, and so a corresponding suspension flow will have a periodic orbit. Thus, to build a flow in a plug without periodic orbits, the transverse holonomy for the flow must be a ``translation'' on an infinite line or region. 

The inspired geometric observation is that such a translation is already available in the Wilson plug, as   the flow on the Reeb cylinder, as pictured in Figure~\ref{fig:reebflow}  below.

\begin{figure}[!htbp]
\centering
 \includegraphics[width=35mm]{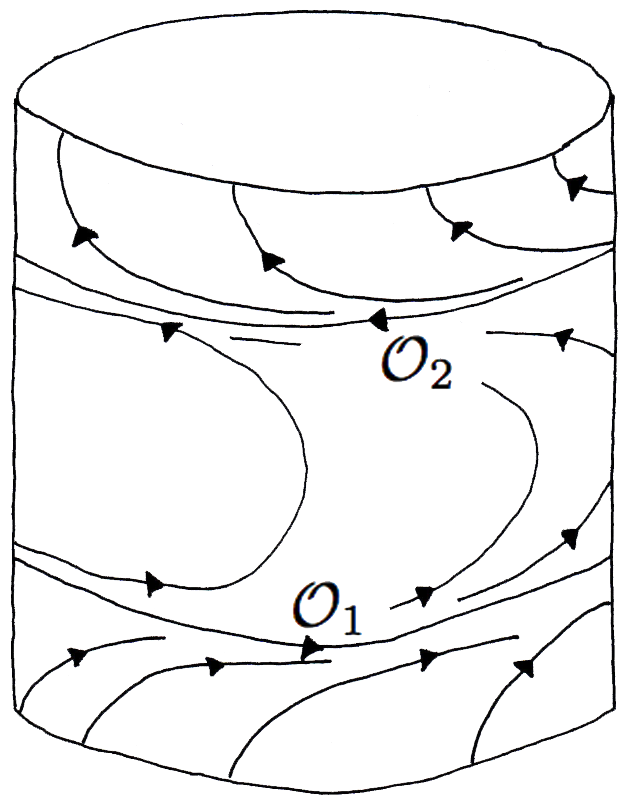}
\caption{Reeb flow on cylinder \label{fig:reebflow}}
\vspace{-6pt}
\end{figure}

The flow of a point below the orbit $\cO_1$ will climb up to the periodic orbit $\cO_1$ and the return map of this flow is a translation on a vertical interval in the cylinder transverse to the flow. A similar dynamical behavior happens for points in the cylinder above the periodic orbit $\cO_2$, but in reverse time.

\begin{figure}[!htbp]
\centering
 {\includegraphics[width=90mm]{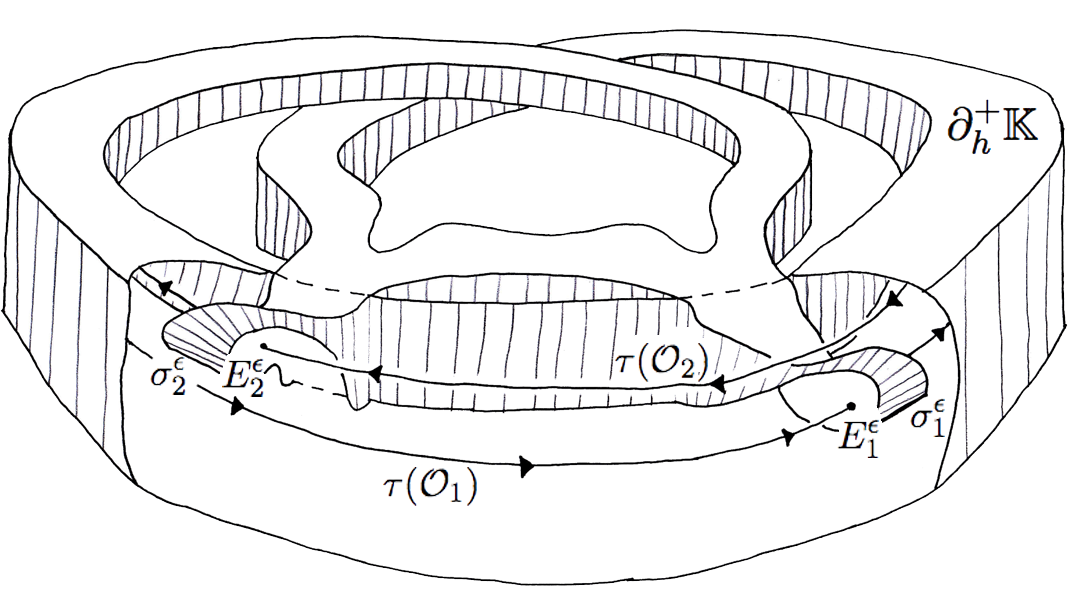}
 \caption{\label{fig:kuperberg} The Kuperberg Plug $\mK_\e$}}
\vspace{-6pt}
\end{figure}

With this point of view, one can then ask: how to get the flow of the periodic orbit to ``break open'' and ``start climbing the holonomy staircase''? That is,   break open the periodic orbit using its own holonomy! 
The solution is the second inspired idea behind Kuperberg's construction,  to insert the flow on the Wilson cylinder  into itself, so that the orbit $\cO_1$ is inserted into one of the flow lines climbing up to $\cO_1$. This is precisely what the insertion maps illustrated in Figure~\ref{fig:kuperberg} below accomplish. To make this work and keep the plug properties, one must first do a ``double twist'' of the embedding of the Reeb cylinder, as illustrated in Figure~\ref{fig:8doblado} below. (See also the illustration on the left side of Figure~\ref{fig:twisted}.)
According to Krystyna, the idea for doing this construction occurred to her during a beer party at the Georgia Topology Conference.

\begin{figure}[!htbp]
\centering
{\includegraphics[width=50mm]{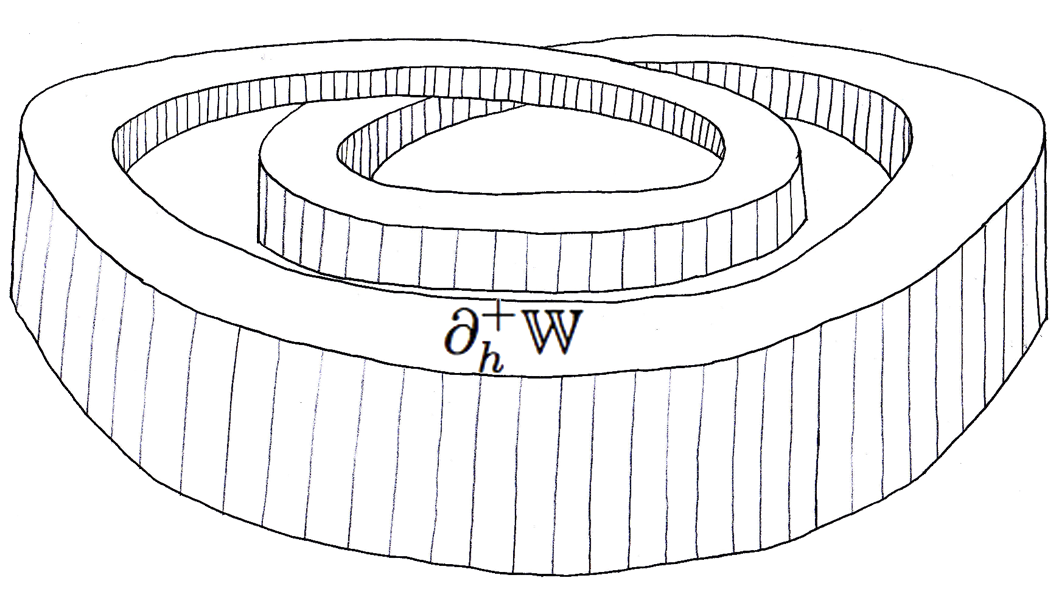}}
\caption{Embedding of Wilson Plug $\mW$ as a {\it folded figure-eight} \label{fig:8doblado} }
\end{figure}

   The essence  of the novel strategy behind the aperiodic property of the flow $\Phi_t$ on a Kuperberg Plug  is perhaps best described by a quote from the paper by Matsumoto \cite{Matsumoto1995}:
\begin{quotation}
  We therefore must demolish the two closed orbits in the Wilson Plug beforehand. But producing a new plug will take us back to the starting line. The idea of Kuperberg is to \emph{let closed orbits demolish themselves}. We set up a trap within enemy lines and watch them settle their dispute while we take no active part. 
\end{quotation}
  The reader is extolled to read the proof that these flows are aperiodic  in \cite{Kuperberg1994}, or in the reports by Ghys \cite{Ghys1995} or Matsumoto \cite{Matsumoto1995}, as the elegance of Kuperberg's idea is revealed in the simplicity of this proof.

  Note that there are many choices of the vector field that can be made for the flow in a Wilson plug. What is essential for Kuperberg's construction is that there are the two closed orbits which attract/repel a set of orbits entering or leaving a face of the plug.

A  Kuperberg plug  can also be constructed for which   its
flow $\cK$ is real analytic. An explicit construction of such a flow
is given in Section~6 of the  paper \cite{Kuperbergs1996} by Greg and Krystyna Kuperberg. There is the added
difficulty that the insertion of the plug in an analytic manifold must
also be analytic, which requires some subtlety. Analytic aperiodic flows are  discussed in   the second author's  Ph.D. Thesis \cite{Rechtman2009}.

 The construction of smooth aperiodic flows on compact manifolds using the Kuperberg method  produces orbits that   never close up, so they wander around a compact region of a plug that has been inserted. There is much that can be said about the dynamical properties of these flows, deduced from the way the flows in the plugs are constructed. 
 These are the \emph{known knowns} about Kuperberg flows, as discussed in Section~\ref{sec-dynamics}.  
 
 There are also many \emph{known unknowns} about Kuperberg flows, which generate interesting open questions about the dynamics of the Kuperberg flows, as discussed in Section~\ref{sec-problems}. There are also the \emph{unknown unknowns}, which are speculations of dynamical phenomena yet to be discovered for this class of smooth flows.

  \section{Aperiodic flow dynamics}\label{sec-dynamics}

 Ghys observed in  \cite{Ghys1995}     that   
a Kuperberg   flow has   zero topological
entropy.  This follows from the remark that   for a smooth flow on a compact 3-manifold,  Katok's   results in   \cite{Katok1980} imply that if a smooth flow has positive entropy, then it must have periodic orbits. Thus, aperiodic flows belong  to a class of smooth dynamical systems which might be considered ``less than chaotic'', or more precisely  ``at worst, slowly chaotic''.

  The flow  in the Wilson plug  in Figure~\ref{fig:Reebcyl} has two periodic orbits, labeled $\cO_1$ and $\cO_2$. Both of these orbits are broken open by the insertion surgery  illustrated in 
 Figure~\ref{fig:kuperberg}, to yield what are called the \emph{special orbits} for the flow. Not as obvious without a closer inspection, is that each of the special orbits is trapped inside of the Kuperberg plug formed by the surgery, and they limit on each other. Thus, their closures form a minimal set for the flow in the plug.  Furthermore, every orbit entering  the plug either escapes from the plug on the opposite face, or else it limits on the closure of a special orbit. Thus, there exists a unique minimal set for the flow in the plug.  If the aperiodic flow is constructed using multiple Kuperberg plugs, then there may be multiple minimal sets, each disjoint from the other. For simplicity of exposition, we assume there is only one plug, and let $\Sigma$ denote the unique minimal set.
 
 Another highly non-obvious result due to  Matsumoto \cite[Proposition~7.2]{Matsumoto1995},   see also the discussion in \cite[Section~8]{Ghys1995}, is that there is an open set of points in the entrance of a Kuperberg plug whose forward orbits limit on the minimal set $\Sigma$. As the minimal set $\Sigma$ is contained in the interior of the plug, this open set consists of non-recurrent  points for the flow. In particular, this implies that there is no invariant probability measure in the Lebesgue measure class, for the   flow in the Kuperberg plug. It is remarked after the proof of \cite[Proposition~18]{Kuperbergs1996} that the Matsumoto argument works for $C^1$-flows, and some but not all piecewise-linear self-insertions  \cite[Section~9]{Kuperbergs1996}.
 
 Remarkably,    Greg Kuperberg showed in Theorem~1 of \cite{KuperbergG1996} that on every closed 3-manifold there is a volume-preserving $C^1$-vector fielld with no closed orbits, obtained by modifying the Schweitzer plug.  He also proved in Theorem~2 of \cite{KuperbergG1996}  that for every closed 3-manifold, there is a volume preserving $C^\infty$-flow with a finite number of periodic orbits. Also in the volume preserving setting, Viktor L. Ginzburg and Ba\c{s}ak G\"urel \cite{GinzburgGurel2003}, constructed a version of Schweitzer's plug that allows to produce examples $C^2$-functions on symplectic 4-manifolds, having a level set whose Hamiltonian vector field has  no periodic orbits. A Hamiltonian vector field always preserves volume and in this case is $C^1$, hence it gives examples of aperiodic volume preserving flows.

 The final general remark about the geometry of the Kuperberg flows is more complicated to explain, but has profound consequences for their dynamics. 
  The Reeb cylinder is an annular region bounded by the two periodic orbits for the Wilson flow,  as pictured on the right side of Figure~\ref{fig:Reebcyl}.
The process of doing flow surgery as illustrated in Figure~\ref{fig:kuperberg} cuts two ``notches'' out of the cylinder. This is the region labeled by $\cR'$ on the left side of  Figure~\ref{fig:notches} below,   called   the  ``notched Reeb cylinder'' as pictured on the right side of Figure~\ref{fig:Reebcyl},  before the cylinder is deformed to appear as in 
 Figure~\ref{fig:8doblado}, which is then twisted  as in  Figure~\ref{fig:kuperberg}.

  \begin{figure}[!htbp]
\centering
\begin{subfigure}[c]{0.4\textwidth}{\includegraphics[height=35mm]{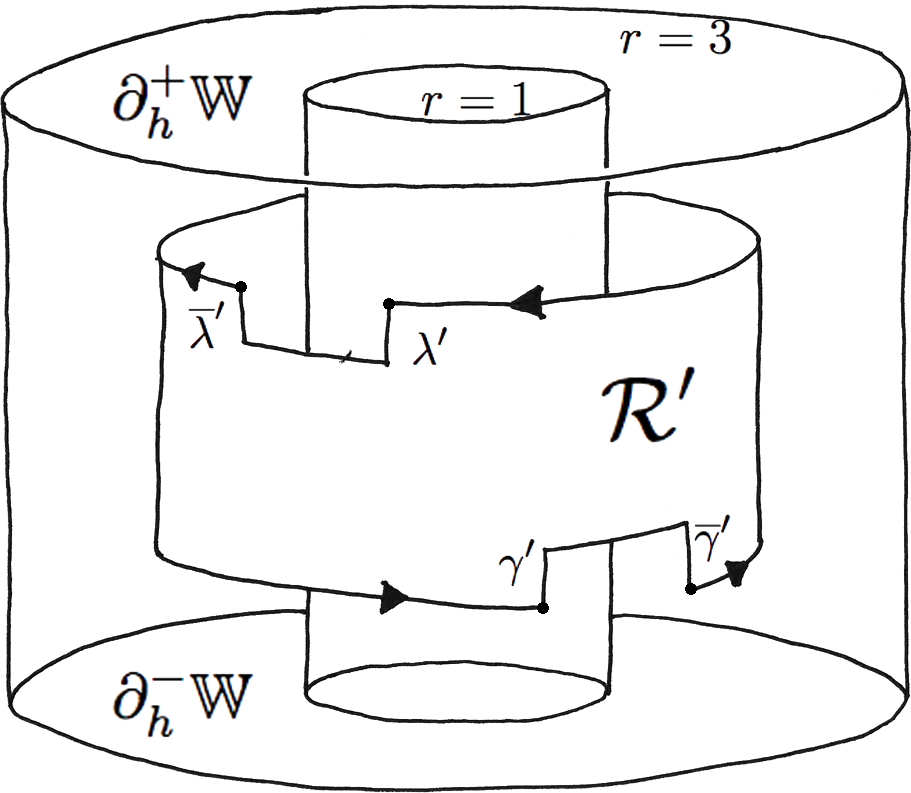}}\end{subfigure}
\begin{subfigure}[c]{0.4\textwidth}{\includegraphics[height=35mm]{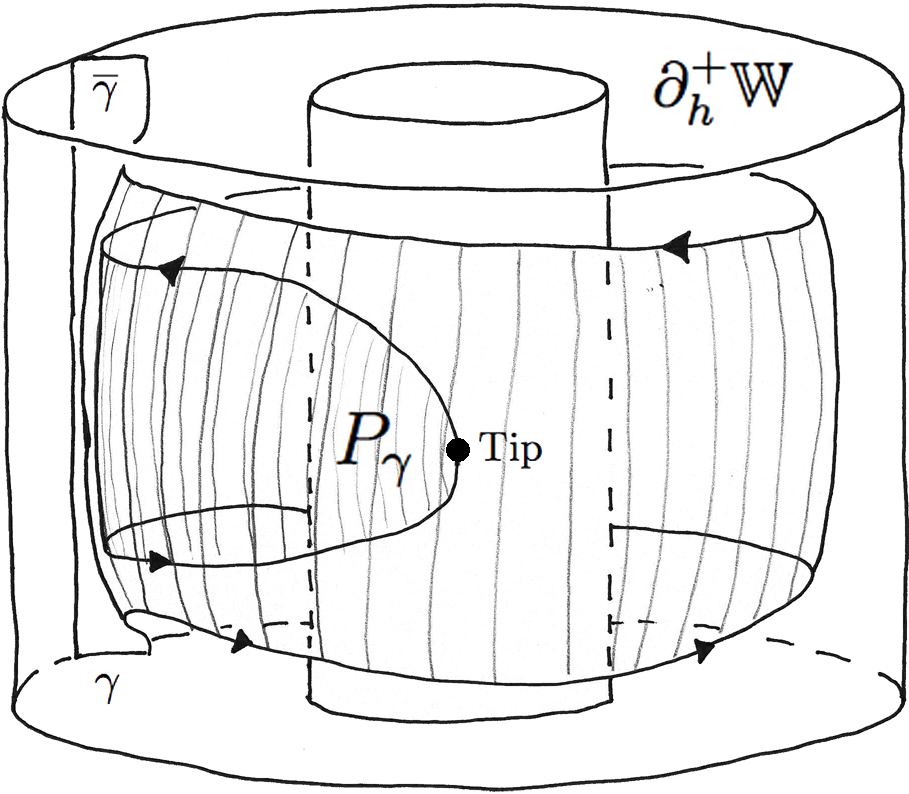}}\end{subfigure}
\caption{\label{fig:notches}   Notch and the flow of its boundary arc}
\vspace{-6pt}
\end{figure}
The Kuperberg flow in the Wilson cylinder preserves the region $\cR'$, called the  \emph{notched Reeb cylinder}, until it reaches a notch, then it continues to flow but the orbits are not longer in the cylinder, as in Figure~\ref{fig:kuperberg}. This produces a tongue-shaped region as illustrated on the right side of  Figure~\ref{fig:notches}. This tongue region continues to flow until it intersects the insertion regions created by the surgery, and the pattern repeats ad infinitum, attaching tongues to the tongues.

\begin{figure}[!htbp]
\centering
{\includegraphics[width=150mm]{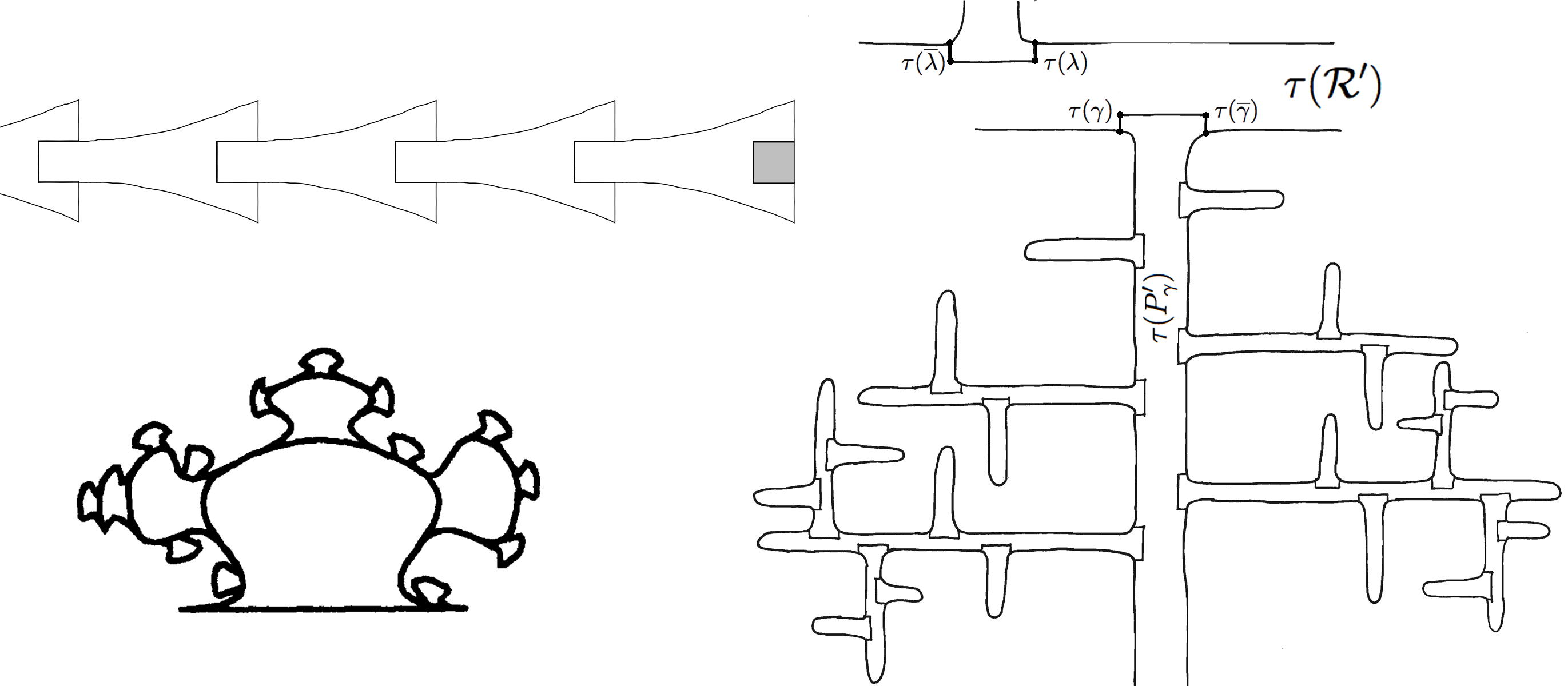}}
\caption{\label{fig:choufleur} The surface $\fM_0$ -   as Serpent,   Chou-fleur, or Propeller}
\vspace{-6pt}
\end{figure}

The complete flow of the \emph{notched Reeb cylinder}  is a surface with boundary embedded in the Kuperberg plug. Figure~\ref{fig:choufleur}   gives three illustrations of this central object for Kuperberg flows, which Siebenmann called a ``chou-fleur'' in his communication \cite{Siebenmann1997}. This was  illustrated by Ghys as a ``fractal-like cluster'' in  \cite{Ghys1995}, as pictured 
in the lower left side of Figure~\ref{fig:choufleur}. In the monograph \cite{HR2016}, the authors 
 called the surface   as laid out on the plane, a  ``propeller'' , which is approximately illustrated by the thickened 2-dimensional, tree-like structure on the right side of Figure~\ref{fig:choufleur}, but again not to scale (all the branches have approximately the same width).  We call this infinite surface $\fM_0$.
     
 The basic observation is that the boundary of $\fM_0$   is the union of  the special orbits for the flow! Thus, to understand the dynamics of the flow, one  analyzes the structure of the  propellers in the Kuperberg plug, which are defined by a recursion process that is generated by the insertion surgery. The   recursion process which defines  the branches of the propeller depends in a sensitive manner on the precise surgery process used to construct the plug. This makes the analysis very delicate, and more will be said about that later.

The Kuperberg flow  preserves the embedded infinite surface   $\fM_0$, whose    boundary       contains the special orbits. The fact that the flow is aperiodic corresponds to the fact that the propellers are infinite. 

The closure $\fM = \overline{\fM_0}$ is a type of 2-dimensional \emph{lamination with boundary} that contains the   minimal set $\Sigma$. Thus, the topological shape of the minimal set $\Sigma$ is dominated by the shape of the lamination $\fM$, which can be estimated using the recursion on the reinsertions of the notched Reeb cylinder for the Kuperberg plug. This  approach to the study of the shape of $\Sigma$ is taken in \cite{HR2016}. Conversely, the inclusion $\Sigma \subset \fM$ suggests possible relations between the shape of $\fM$ and the dynamical properties of the flow, as discussed in Section~\ref{sec-problems}.

We describe the technical setting for comparing $\Sigma$ and $\fM$.
The drawing on the left side of Figure~\ref{fig:twisted}  below gives an expanded view of the  insertions illustrated in Figure~\ref{fig:kuperberg}.  A key point is that the image of a boundary orbit, denoted by $\tau(\cO_1)$ in the illustration, has image under the insertion as a twisted curve that is tangent to the cylinder at the boundary curve $\cO_1$, but intersects transversally $\tau(\cO_1)$.  This results in a picture as on the right side of Figure~\ref{fig:twisted}, where the two thickened curves are the boundaries of the cylinder and its insertion. 

The graph on the right side of Figure~\ref{fig:twisted}  is called the \emph{radius function} (actually, what is illustrated is the inverse of this function) and   the ``quadratic shape'' of the graph is the basis for showing that the Kuperberg flow is aperiodic. Understanding these remarks takes some effort, and is explained in the papers \cite{Kuperberg1994,Ghys1995,Matsumoto1995,HR2016}, but reveals the great beauty of the Kuperberg construction.

  The authors introduced the notion of a \emph{generic} Kuperberg flow in the work \cite{HR2016}, which formulates a collection of optimal conditions on the choices for the construction of the flow. One of these conditions is that  the radius function for the insertion is actually ``quadratic'', as pictured on the right side of Figure~\ref{fig:twisted}.
  We   show in \cite{HR2016} that for a generic Kuperberg flow, there is equality $\Sigma = \fM$, and thus one can study the shape of the minimal set for the flow by studying the properties of the 2-dimensional lamination with boundary $\fM$.
  
 \begin{figure}[!htbp]
\centering
{\includegraphics[width=44mm]{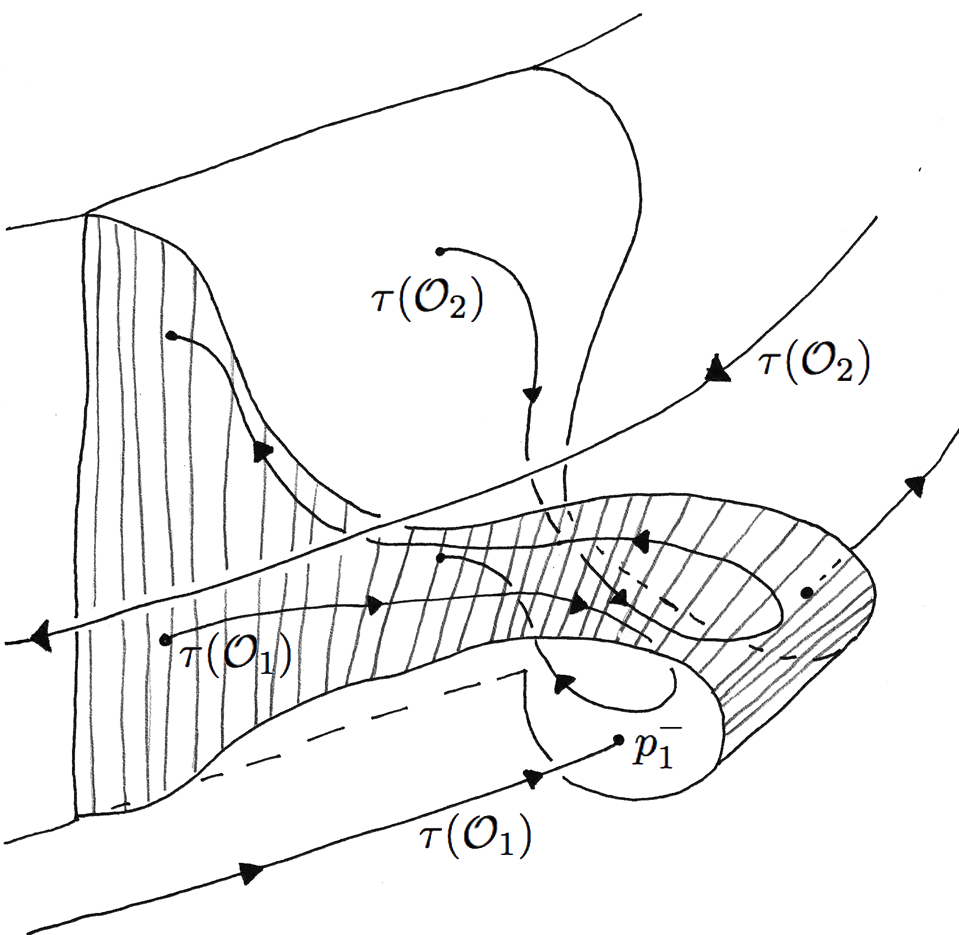}\hspace{20pt} \includegraphics[width=50mm]{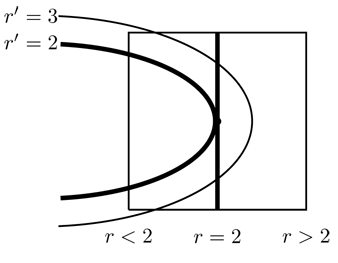}}
\caption{The image of $L_1\times [-2,2]$ under $\sigma_1$ and the radius function\label{fig:twisted} }
\vspace{-6pt}
\end{figure}
 
 \eject
 
  \section{Non-generic flows}\label{sec-problems}

Ghys wrote in his  1995 survey \cite[page 302]{Ghys1995}: 
\begin{quotation}
\emph{Par ailleurs, on peut construire beaucoup de pi\`eges de Kuperberg et il n'est pas clair qu'ils aient la m\^eme dynamique.}
\end{quotation}

In this section, we discuss   some  ways in which the dynamical properties of  an aperiodic smooth Kuperberg plug   vary with  the choices made in its construction. The properties considered include the \emph{shape} and the \emph{Hausdorff dimension} of the minimal set $\Sigma$,  and the  \emph{slow entropy} of the flow.

The construction of a Kuperberg flow begins with the choice of  a flow $\cW$ on a cylinder,  with an invariant cylinder bounded by the two periodic orbits that is called    the \emph{Reeb cylinder}, as illustrated on the left side of  Figure~\ref{fig:Reebcyl}. 
 This flow is then extended smoothly to a thickened cylinder, or interval times an annulus, as on the right side of  Figure~\ref{fig:Reebcyl}. There are further symmetry conditions imposed on the flow $\cW$, as given in \cite[Section~3]{Kuperberg1994}, or for  example as specified by (P1) to (P4) in \cite[Section~2]{HR2016}. The symmetry conditions still  allow for a wide variation of choices for the  flow  $\cW$.
  
 The flow of the Kuperberg plug, that  results from the insertion process, consists of  orbit segments of the gradient like flow $\cW$ which are ``patched together''. This view of the orbits, as a union of finite flows, is key to the analysis of the Kuperberg flow in the works of Ghys \cite{Ghys1995} and Matsumoto \cite{Matsumoto1995}. For a different take on this, 
 Section~5 of the   paper \cite{Kuperbergs1996} describes the result of the self-insertion patching, as the analog of a ``stacking subroutine'' for a computer program.  Then the fact that no orbit is periodic corresponds to showing that the routine does not terminate for the given starting point. So in essence, the solution to the Seifert Conjecture is realized by building a ``dynamical computer program'' that is non-terminating.
  
 This  patching of flow segments is occurring orbitwise, so if there is to be any hope of organizing this process in a unified description for the entire flow on the plug, one needs to impose assumptions on the flow $\cW$ that its orbits are as uniform as possible.
 Such additional conditions on $\cW$ are formulated in the works \cite{Kuperberg1994,Kuperbergs1996}, and even more restrictive conditions are imposed in \cite{HR2016}. 
   
 The second step in the construction of a Kuperberg flow is to  make a choice for each of the two insertion maps, as pictured in Figure~\ref{fig:kuperberg} and Figure~\ref{fig:twisted}.  In this step,   the embedding of the Reeb cylinder makes a transverse contact at the periodic orbit, either the bottom or the top orbit, so as to break open these closed orbits. The result is that the special orbits which result by breaking open the closed orbits, they intersect the   vertex of the parabola pictured on the right side of Figure~\ref{fig:twisted}. They  return infinitely often  to a neighborhood of the vertex of the parabola as a result of the translation holonomy of the flow. Hence, the   germinal shape at its vertex  of the parabolic map     has a strong influence on the dynamics of the special orbit, and so on the dynamics of the entire flow.
  
 A basic question is,  when does the minimal set satisfy $\Sigma = \fM$? Ghys gave a construction of a flow on a plug  and sketched the proof that $\Sigma = \fM$ for this flow, at the end of Section~7 in \cite{Ghys1995}. The Kuperbergs gave an explicit real analytic flow on a plug for which $\Sigma = \fM$   in their joint paper \cite{Kuperbergs1996}. In both cases, the idea of the proof is to show that the closure of a special orbit contains the notched Reeb cylinder $\cR'$, hence contains the flow of $\fM_0$ and thus equals its closure $\fM$.
  
  The authors began our investigations of the shape of the minimal set in the Kuperberg flows in 2010, and for several years after that,  every time we met Krystyna we had to admit we couldn't answer one of the mysteries of the Kuperberg construction: 
 \begin{prob}\label{prob-denjoy}
 Is the minimal set \emph{always} two dimensional for a smooth Kuperberg flow?
 \end{prob}
 This is the conclusion for the examples by Ghys in \cite{Ghys1995} and in the  joint paper \cite{Kuperbergs1996}, but the situation in general was not known. 
 In order to answer Krystyna's question, we   introduced the notion of a ``generic Kuperberg flow'' in  \cite{HR2016}, which assumes that the flow satisfies the conditions  of Hypothesis~12.2 and Hypothesis~17.2 formulated in that work.
  Hypothesis~12.2  assumes in particular  that  the vertical component of the vector field $\cW$ on the Reeb cylinder has quadratic germ near the two periodic orbits at which its ``vertical component'' must vanish. This condition is used to show that the closure of the special orbits contain $\cR'$. But it is unknown, for example, that if the vector field vanishes to higher order at the periodic orbits, can one still show that $\cR'$ is contained in the closure of a special orbit?

 The equality  $\Sigma = \fM$ can be considered as a type of \emph{Denjoy Theorem} for 2-dimensional laminations. The usual 
Denjoy Theorem states  that the closure of an orbit for a   $C^2$-flow without closed orbits on the 2-torus is all of  $\mT^2$. Then the equality $\Sigma = \fM$ can be considered   as analogous, for  it states that the closure of an orbit is not a 1-dimensional submanifold of $\fM$, but fills up its  2-dimensional leaves. The known proofs of this conclusion cited above, they  all make  assumptions on both the flow $\cW$ and on the insertion map. Formulating  general criteria which suffice to imply the equality $\Sigma = \fM$ has proven difficult.

 If  $\Sigma \subset \fM$ is a proper inclusion, then the topological type of $\Sigma$ and $\fM$ will differ, for example they could have distinct shapes. A $C^1$-flow on $\mT^2$ without closed orbits provides a good model for this difference. If the minimal set $\fD$ for the flow is not all of $\mT^2$ then $\fD$ is a 1-dimensional continuum, as used in the construction of the Schweitzer plug. If $p \in \mT^2 - \fD$  then there is a retract of the open punctured torus $\mT^2 - \{p\}$ onto $\fD$, so that $\fD$ has the shape of a wedge of two circles. It seems a very interesting problem to compare the shapes of the  closed invariant sets  for a general aperiodic flow:
   \begin{prob}\label{prob-punctured}
 Suppose the Kuperberg flow in a plug  has minimal set $\Sigma \subset \fM$ which is a proper inclusion. 
 What is the relation between the shapes of $\Sigma$ and   $\fM$?
 \end{prob}
Theorem~19 in \cite{Kuperbergs1996} gives the construction of aperiodic PL plugs for which $\Sigma$ is 1-dimensional. In addition, there is a discussion of the symbolic dynamics for these special flows they construct. 
The consideration of the topological type of $\Sigma$, in addition to the known unknowns about its shape, one senses that its  study    leads  into the realm of the unknown unknowns, that new dynamical phenomena will be discovered.
 
 We next consider another dynamical aspect of Kuperberg flows, their entropy and the Hausdorff dimension of their minimal sets. While Ghys observed in \cite{Ghys1995} that an aperiodic flow must have entropy equal to zero, using a well-known deep result of Katok \cite{Katok1980}, none-the-less the dynamical behavior of the flow appears to be chaotic upon closer inspection.
 
 The geometry of the propeller as pictured in Figure~\ref{fig:choufleur}, suggests that the flow in the minimal set $\Sigma$ should exhibit exponential behavior, as it follows the boundary of the propeller which appears to have exponentially growing area. This intuition is flawed, though, as the flow cannot simply travel down the core of the propeller, but must instead travel on the boundary to get to the extremities of the tree. The time required to get to the extremes of the tree is not linear, but grows approximately as the square of the distance to be traveled. This observation is the basis for Theorem~21.10 in \cite{HR2016}: 
 \begin{thm}\label{thm-slow}
The ``slow entropy''   of a generic Kuperberg flow is positive for exponent $\alpha = 1/2$. 
 \end{thm}
 The notion of slow entropy was introduced by   Katok and Thouvenot  \cite{KatokThouvenot1997}. The exponent $\alpha=1/2$ is essentially saying that the entropy would be positive if the time variable is sped up by taking the square of time. We also note that for this result, the generic hypotheses include an extra condition, the assumption that the insertion maps for the construction of the plug have ``slow growth'' themselves. This is a technical property, but its requirement  emphasizes that the proof of Theorem~\ref{thm-slow} is quite technically involved. In particular, the proof uses an estimate on the growth rates of leaves in the lamination $\fM$, that they grow at subexponential rate $1/2$ as well.
Define the growth rate of $\fM$ as the exponential growth rate of the dense leaf $\fM_0$, so for a generic flow this rate is $1/2$. Also, define the entropy dimension $HD(\fM)$ of $\fM$ as the least upper bound of the exponents $0 \leq \alpha \leq 1$ such that the $\alpha$-slow entropy of the flow is positive \cite{deCarvalho1997}. Thus, for a generic flow we have $HD(\fM) \geq 1/2$.
There are many questions one can ask about these  invariants, for example:

    \begin{prob}\label{prob-entropy}
What  range of values for $HD(\fM)$   can be realized by aperiodic flows on 3-manifolds? Is there an aperiodic flow with $HD(\fM) = 0$?
\end{prob}
 
 For a Kuperberg flow, one can also ask about the  Hausdorff dimension of  the lamination $\fM$.
This problem seems almost intractable, but Daniel Ingebretson developed an approach to calculating this dimension for a particular class of flows in his thesis work \cite{Ingebretson2018a,Ingebretson2018b}. 

 \begin{thm}\label{thm-dim}\cite{Ingebretson2018a}
The Hausdorff dimension   of the minimal set for  generic Kuperberg flow satisfies $2 < d_h(\fM) < 3$.
 \end{thm}
The   dimension of $\fM$ must be at least $2$, as $\fM$ is a union of $2$-dimensional ``leaves'' by the generic assumption. But the fact that the dimension is greater than 2 is a measure of the transverse complexity of $\fM$. The actual Hausdorff dimension appears to depend in a   sensitive manner on the choices made. The thesis work also gives a method to numerically calculate this number.
We can then formulate a question which falls into the category of an unknown unknown about Kuperberg flows:
    \begin{prob}\label{prob-coupling}
Is there a relation between the slow entropy of a Kuperberg flow, and the Hausdorff dimension of its unique minimal set? Is there a relation if the flow is assumed to be generic?
\end{prob}

We should also mention a question that belongs to the class of ``unknown unknown'' problems. 
  
  \begin{prob}
How do the dynamics of a real analytic Kuperberg flow in a plug differ from the general smooth case?
Are there further restrictions on the shape of the minimal sets in the analytic case?
\end{prob}

Another remarkable fact about the Kuperberg flows, is that they   lie at the ``boundary of chaos'' in the $C^{\infty}$-topology on flows.   The idea behind this remark is that when making the insertion as in Figure~\ref{fig:twisted}, one can stop the insertion ``too soon''. That is, if the two cylinders do not make contact, then the insertion does not break open the periodic orbits for the Wilson flow $\cW$. We show in \cite{HR2018a} that the flow  for such a truncated insertion again has two periodic orbits. However, the lengths of these orbits tends to infinity as the insertion limits to one which breaks open the orbits. 
Thus, we obtain  smooth  families  of variations of the
  Kuperberg plug with simple dynamics,  and exactly two periodic orbits,
  and the limit of the family is an aperiodic flow.  
  
  In such a family,
  the period of the periodic orbits must blow up at the limit.
 Palis and Shub in  \cite{PalisPugh1975} asked whether this dynamical phenomenon can occur in families of smooth flows on closed manifolds, and called a closed orbit whose length ``blows up'' to infinity under deformation 
a ``blue sky catastrophe''.  
The first examples of a family of flows with this property was found by   Medvedev in  \cite{Medvedev1980}. The constructions in \cite{HR2018a} show that deformations of Kuperberg flows provide a new class of examples.  The work of  A.~Shilnikov, L.~Shilnikov and D.~Turaev in  \cite{SST2005} further discusses ``blue sky catastrophe'' phenomenon.

 On the other hand,  when making the insertion as in Figure~\ref{fig:twisted}, one can take the insertion ``too far''. That is, the two cylinders  intersect not along one arc, {\it but along two arcs}.  We show in \cite{HR2018a} that the flow  for these over-extended insertions have countably many embedded horseshoes, so an abundance of hyperbolic behavior. Moreover, as the insertion maps are deformed to one that just makes contact at a single point, then these horseshoes increase in number. The dynamics of the generic Kuperberg flow is thus  the limit of the dynamical chaos in a continuous family of horseshoes for neighboring smooth flows.

 The immediate conclusion is that the Kuperberg flows are not stable in the $C^{\infty}$-topology on flows.

 The recitation of remarkable properties of the class of flows created by Kuperberg could continue, as they are zero entropy dynamical systems, which are simply not boring \cite{HR2018b}! Instead, we want to also discuss a variety of other constructions and results following from Seifert's Conjecture, and the questions they engender.

  \section{Decorated flows}\label{sec-speculations}
  
 As explained by Ghys \cite{Ghys1995}, the construction of non-singular aperiodic flows is particularly complicated in dimension 3. The Seifert conjecture, that is false in full generality, can be decorated with invariant structures: we can ask, if the flow preserves something, then must it have a periodic orbit? On one hand, there are families of non-singular flows that always have periodic orbits, on the other hand we have the examples discussed previously in this text. 
  Here we discuss this broader subject of flows on closed 3-manifolds, we review some results on the existence of periodic orbits.  We do not intend to make an exhaustive list.
  
  Let us start with  the class of Reeb vector fields associated to a contact form. It has been shown by   Hofer in \cite{Hofer1993}, that the flow of a Reeb vector field on $\mS^3$  must have a   periodic orbit, and this was generalized to every closed 3-manifold by    Taubes  in \cite{Taubes2007}. But we know more, in fact every Reeb vector field on a closed 3-manifold has at least two periodic orbits as proved by Cristofaro-Gardiner and Hutchings \cite{Hut-CG2016} and the examples with exactly two periodic orbits are completely understood \cite{TWO2021}.  The second author,   together with Colin and Dehornoy, proved that if a Reeb field is non-degenerate then it has either two or infinitely many periodic orbits \cite{CDR2020}. 
  
A Reeb vector field always preserves a volume form. A big open question in the subject, is whether a $C^\infty$ volume preserving flow on a 3-manifold has a periodic orbit. It seems impossible to mimic Krystyna's construction of a plug in this setting, but we can do a Wilson plug for volume preserving flows. The lack of examples of possible minimal sets, makes it also impossible to mimic other plug constructions.

In the case of Reeb flows, Alves and Pirnapasov \cite{AlvesPir2021} proved recently the existence of some knots  such that if a Reeb flow has a periodic orbit realizing the knot then it has positive entropy. In particular, by the already mentioned theorem of Katok \cite{Katok1980}, such Reeb flows have lots of periodic orbits. It is an observation that such a result cannot hold for categories of flows for which we can construct plugs. Using a suitable Wilson plug one can prove the following:
\bigskip

\begin{thm} Let $K$ be a given knot. Every 3-manifold that admits a zero entropy flow without fixed points, admits a zero entropy flow without fixed points having (at least) two periodic orbits whose knot type is $K$.
Moreover, if the original flow is generated by a Reeb vector field, then the new flow  preserves a plane field that fails to be a contact structure along the two periodic orbits of knot type $K$.
\end{thm}

 The existence of periodic orbits was extended to geodesible volume preserving flows (also known as Reeb vector fields of stable Hamiltonian structures) on manifolds that are not torus bundles over the circle by Hutchings and Taubes \cite{HutchingsTaubes2009,Hutchings2010}, and Rechtman \cite{Rechtman2010}, and for real analytic geodesible flows by Rechtman \cite{Rechtman2010}. The existence of a periodic orbit is also established for real analytic solutions of the Euler equation by Etnyre and Ghrist  \cite{EtnyreGhrist2000}.  Geodesible volume preserving flows are solutions to the Euler equation, but we do not know if the last ones have periodic orbits. We do know that it is impossible to construct plugs whose flow satisfies the Euler equation \cite{PSRTL2021}.


\begin{thebibliography}{10}

\bibitem{AlvesPir2021}
{M.~Alves and A.~Pirnapasov},
\newblock{\it Reeb orbits that force topological entropy}, Ergodic Theory and Dynamical Systems, 2021. 

\bibitem{Borsuk1968}
{K.~Borsuk},
\newblock {\it Concerning homotopy properties of compacta},
\newblock {\bf  Fund. Math.}, 62:223--254, 1968.


\bibitem{Borsuk1975}
{K.~Borsuk},
\newblock {\bf Theory of shape},
\newblock Monografie Mat., vol. 59, Polish Science Publ., Warszawa, 1975..


\bibitem{CDR2020} 
V. Colin, P. Dehornoy and A. Rechtman, {\it On the existence of supporting broken book decompositions for contact forms in dimension $3$},  {arXiv:2001.01448}.

\bibitem{TWO2021}
{D.~Cristofaro-Gardiner, U.~Hryniewicz, M.~Hutchings and H.~Liu},
\newblock{\it Contact three-manifolds with exactly two simple Reeb orbits}, {arXiv:2102.04970}.

\bibitem{Hut-CG2016}
{D.~Cristofaro-Gardiner and M.~Hutchings}, 
\newblock{\it From one Reeb orbit to two},
\newblock {\bf J. Differential Geom.} 102, no. 1, 25--36, 2016.

\bibitem{deCarvalho1997}
{A.~de Carvalho},
\newblock {\it Entropy dimension of dynamical systems},
\newblock {\bf Portugal. Math.}, 54:19--40, 1997.

   
 \bibitem{Denjoy1932}
{A.~Denjoy},
\newblock {\it Sur les courbes d\'efinies par les \'equations diff\'erentielles \`a la surface du tore,},
\newblock {\bf J. Math. Pure Appl.}, 11:333-375, 1932.


\bibitem{EtnyreGhrist2000}
{J.~Etnyre and R.~Ghrist},
\newblock {\it Contact topology and hydrodynamics. {I}. {B}eltrami fields and the {S}eifert conjecture},
\newblock {\bf  Nonlinearity}, 13:441--458, 2000.

 

\bibitem{Ghys1995}
{{\'E}~Ghys},
\newblock {\it Construction de champs de vecteurs sans orbite p\'eriodique (d'apr\`es {K}rystyna {K}uperberg)},
\newblock {S{\'e}minaire Bourbaki, Vol. 1993/94, Exp.\ No.\ 785},
\newblock {\bf Ast\'erisque}, Vol. 227: 283--307, 1995.



         
\bibitem{GinzburgGurel2003}
{V.L.~Ginzburg and B.Z.~G\"{u}rel},
\newblock {\it A {$C^2$}-smooth counterexample to the {H}amiltonian {S}eifert conjecture in {$\Bbb R^4$}},
\newblock {\bf Ann. of Math. (2)}, 158:953--976, 2003.

         
\bibitem{Handel1980}
{M.~Handel},
\newblock {\it One-dimensional minimal sets and the {S}eifert conjecture},
\newblock {\bf Ann. of Math. (2)}, 111:35--66, 1980.



\bibitem{Harrison1988}
{J.C.~Harrison},
\newblock {\it {$C^2$} counterexamples to the {S}eifert conjecture},
\newblock {\bf Topology}, 27:249--278, 1988.


\bibitem{Hofer1993}
{H.~Hofer},
\newblock {\it Pseudoholomorphic curves in symplectizations with applications   to the {W}einstein conjecture in dimension three},
\newblock {\bf Invent. Math.}, 114:515--563, 1993.



\bibitem{HR2016}
{S.~Hurder and A.~Rechtman},
\newblock {\it The dynamics of  generic Kuperberg flows},
\newblock  {\bf Ast\'erisque}, Vol. 377:1--250, 2016.

\bibitem{HR2018a}
{S.~Hurder  and A.~Rechtman},
\newblock {\it Aperiodicity at the boundary of chaos},
\newblock  {\bf Ergodic Theory Dyn. Systems}, 38:2683--2728, 2018.

 

\bibitem{HR2018b}
{S.~Hurder  and A.~Rechtman},
\newblock {\it Perspectives on {K}uperberg flows},
\newblock  {\bf Topology Proc.},   51:197--244, 2018.

 
\bibitem{HutchingsTaubes2009}
{M.~Hutchings and C.~Taubes},
\newblock {\it The {W}einstein conjecture for stable {H}amiltonian structures},
\newblock {\bf Geom. Topol.}, 13:901--941, 2009.

      
\bibitem{Hutchings2010}
{M.~Hutchings},
\newblock {\it Taubes's proof of the {W}einstein conjecture in dimension three},
\newblock {\bf Bull. Amer. Math. Soc. (N.S.)}, 47:73--125, 2010.
    
\bibitem{Ingebretson2018a}
{D.~Ingebretson},
\newblock {\it Hausdorff dimension of generic Kuperberg minimal sets},
\newblock {\bf Thesis}, University of Illinois at Chicago, 2018.
    
\bibitem{Ingebretson2018b}
{D.~Ingebretson},
\newblock {\it Hausdorff dimension of Kuperberg minimal sets},
\newblock {preprint}, {arXiv:1801.04034}.
    
     
\bibitem{Katok1980}
{A.~Katok},
\newblock {\it Lyapunov exponents, entropy and periodic orbits for diffeomorphisms},
\newblock {\bf Inst. Hautes \'Etudes Sci. Publ. Math.}, 51:137--173, 1980.

 
\bibitem{KatokThouvenot1997}
{A.~Katok and J.-P.~Thouvenot},
\newblock {\it Slow entropy type invariants and smooth realization of  commuting measure-preserving transformations},
\newblock {\bf Ann. Inst. H. Poincar\'e Probab. Statist.}, 33:323--338, 1997.



\bibitem{KuperbergReed1981}
{K.~Kuperberg and C.~Reed},
\newblock {\it A rest point free dynamical system on {${\mathbb R}^{3}$} with uniformly bounded trajectories},
\newblock {\bf Fund. Math.}, 114:229--234, 1981.

\bibitem{KuperbergReed1989}
{K.~Kuperberg and C.~Reed},
\newblock {\it A dynamical system on {${\mathbb R}^3$} with uniformly bounded  trajectories and no compact trajectories},
\newblock {\bf Proc. Amer. Math. Soc.}, 106:1095--1097, 1989.

  
\bibitem{KKPR1993}
{K.~Kuperberg, W.~Kuperberg, P.~Minc and C.~Reed},
\newblock {\it Examples related to {U}lam's fixed point problem},
\newblock {\bf Topol. Methods Nonlinear Anal.}, 1:173--181, 1993.

  

\bibitem{Kuperberg1994}
{K.~Kuperberg},
\newblock {\it A smooth counterexample to the {S}eifert conjecture},
\newblock {\bf Ann. of Math. (2)}, 140:723--732, 1994.

\bibitem{Kuperbergs1996}
{G.~Kuperberg and K.~Kuperberg},
\newblock {\it Generalized counterexamples to the {S}eifert conjecture},
\newblock {\bf Ann. of Math. (2)}, 144:239--268, 1996.

\bibitem{KuperbergG1996}
{G.~Kuperberg},
\newblock {\it A volume-preserving counterexample to the {S}eifert conjecture},
\newblock {\bf Comment. Math. Helv.}, 71:70--97, 1996.

 

\bibitem{Kuperberg1998}
{K.~Kuperberg},
\newblock {\it Counterexamples to the {S}eifert conjecture},
\newblock in {\bf Proceedings of the {I}nternational {C}ongress of
              {M}athematicians, {V}ol. {II} ({B}erlin, 1998)}, 
\newblock {Doc. Math.},{Extra Vol. II}, 1998, pages {831--840}.




\bibitem{Kuperberg1999}
{K.~Kuperberg},
\newblock {\it Aperiodic dynamical systems},
\newblock {\bf Notices Amer. Math.  Soc.}, 46, no. 9:1035--1040, 1999.


\bibitem{Kuperberg2021}
{K.~Kuperberg},
\newblock {\it personal communication}, November 15, 2021.


   
\bibitem{MardesicSegal1982}
{S.~Marde{\v{s}}i{\'c} and J.~Segal},
\newblock {\bf Shape theory: The inverse system approach},
\newblock  {North-Holland Math. Library, Vol. 26},  North-Holland Publishing Co., Amsterdam, 1982.


\bibitem{Mardesic1999}
{S.~Marde{\v{s}}i{\'c}},
\newblock {\it Absolute neighborhood retracts and shape theory},
\newblock in {\bf History of {T}opology},  pages 241--269, North-Holland Publishing Co., Amsterdam, 1999.


\bibitem{Matsumoto1995}
{S.~Matsumoto},
\newblock {\it {K}uperberg's {$C^\infty$} counterexample to the {S}eifert conjecture},
\newblock {\bf S\=ugaku}, Mathematical Society of Japan, 47:38--45, 1995.
\newblock {Translation:} {\bf Sugaku Expositions}, 11:39--49, 1998, 
\newblock {Amer. Math. Soc.} 
 
 
   
\bibitem{Mauldin2015}
{R.D.~Mauldin},
\newblock {\bf Mathematics from the Scottish Caf\'{e} with selected problems from
              the new Scottish Book,
              Including selected papers presented at the Scottish Book
              Conference held at North Texas University, Denton, TX, May
              1979},
\newblock  {Birkh\"{a}user/Springer, Cham},  2015.

 
\bibitem{Medvedev1980}
{V.S.~Medvedev},
\newblock {\it A new type of bifurcations on manifolds},
\newblock {\bf Mat. Sb. (N.S.)}, 113(155):487--492, 496, 1980.


 \bibitem{PalisPugh1975}
{J.~Palis and C.~Pugh},
\newblock {\it Fifty problems in dynamical systems},
\newblock  in {\bf Dynamical systems---{W}arwick 1974 ({P}roc. {S}ympos. {A}ppl.
              {T}opology and {D}ynamical {S}ystems, {U}niv. {W}arwick,
              {C}oventry, 1973/1974; presented to {E}. {C}. {Z}eeman on his
              fiftieth birthday)},  
\newblock {Lecture Notes in Math., Vol. 468}, pages 345--353, Springer, Berlin, 1975

\bibitem{PSRTL2021}
{D.~Peralta Salas, A.~Rechtman and F.~Torres de Lizaur},
\newblock{\it A characterization of 3{D} steady Euler flows using commuting zero-flux homologies},
\newblock {\bf Ergodic Theory Dynam. Systems} 41 , no. 7, 2166--2181, 2021.

\bibitem{PercelWilson1977}
{P.B.~Percell and F.W.~Wilson, Jr.},
\newblock {\it Plugging flows},
\newblock {\bf Trans. Amer. Math. Soc.}, 233:93--103, 1977.


 \bibitem{Rechtman2009}
{A.~Rechtman},
\newblock {\it Pi\`{e}ges dans la th\'{e}orie des feuilletages:  exemples et contre-exemples},
\newblock  {\bf Th\`{e}se}, l'Universit\'{e} de Lyon -- \'{E}cole Normale Sup\'{e}rieure de Lyon, 2009.


 \bibitem{Rechtman2010}
{A.~Rechtman},
\newblock {\it Existence of periodic orbits for geodesible vector fields on  closed 3-manifolds},
\newblock  {\bf Ergodic Theory Dynam. Systems}, 30:1817--1841, 2010.

       
       
\bibitem{Schweitzer1974}
{P.A.~Schweitzer},
\newblock {\it Counterexamples to the {S}eifert conjecture and opening closed leaves of foliations},
\newblock  {\bf Ann. of Math. (2)}, 100:386--400, 1974.

\bibitem{Schweitzer1975}
{P.A.~Schweitzer},
\newblock {\it Compact leaves of foliations},
\newblock  in {\bf Proceedings of the {I}nternational {C}ongress of
              {M}athematicians ({V}ancouver, {B}. {C}., 1974), {V}ol. 1}, {543--546}, 1975.
      
      
\bibitem{Seifert1950}
{H.~Seifert},
\newblock {\it Closed integral curves in {$3$}-space and isotopic two-dimensional deformations},
\newblock {\bf Proc. Amer. Math. Soc.}, 1:287--302, 1950.

\bibitem{SST2005}
{A.~Shilnikov, L.~Shilnikov and D.~Turaev},
\newblock {\it Blue-sky catastrophe in singularly perturbed systems},
\newblock {\bf Mosc. Math. J.}, 5:269--282, 2005.

 
\bibitem{Siebenmann1997}
{L.~Siebenmann},
\newblock {\it Le paradigm du serpent: explique les contre-exemples de {K.} {Kuperberg} \`{a} la conjecture de {Seifert}},
\newblock {manuscript priv\'{e}e}, d\'{e}cembre, 1997.
 


 \bibitem{Taubes2007}
{C.~Taubes},
\newblock {\it The {S}eiberg-{W}itten equations and the {W}einstein conjecture},
\newblock {\bf Geom. Topol.}, 11:2117--2202, 2007.

     
\bibitem{Wilson1966}
{F.W.~Wilson, Jr.},
\newblock {\it On the minimal sets of non-singular vector fields},
\newblock {\bf Ann. of Math. (2)}, 84:529--536, 1966.

\end{thebibliography}
\end{document}